\documentclass[oneside,10pt]{amsart}
\usepackage{amsmath,amsfonts, latexsym,amssymb}
\usepackage[active]{srcltx}

\theoremstyle{definition}

\numberwithin{equation}{section}


\begin{document}

\title[A characterization of maximal ideals in  the Fr\'{e}chet algebras
 \ldots  ]{A characterization of maximal ideals in  the Fr\'{e}chet algebras  
of holomorphic functions $F^p$ $(1<p<\infty$)}

\author{Romeo Me\v strovi\' c}

\address{University of Montenegro, Maritime Faculty Kotor,
\,\, Dobrota 36, 85330 Kotor, Montenegro, 
e-mail: {\tt romeo@ucg.ac.me}}

\begin{abstract}
The space $F^p$ ($1<p<\infty$)  consists of  
 all   holomorphic  functions $f$  on the open unit disk $\Bbb D$ for which    
   $ \lim_{r\to 1}(1-r)^{1/q}\log^+M_{\infty}(r,f)=0,$ 
where $M_{\infty}(r,f)=\max_{\vert z\vert\le r}\vert f(z)\vert$ with $0<r<1$.
 Stoll  \cite[Theorem 3.2]{s} proved that the space $F^p$ 
with the topology given by the  family of seminorms 
$\left\{\Vert \cdot\Vert_{q,c}\right\}_{c>0}$ defined for 
$f\in F^q$ as 
  $\Vert f\Vert_{q,c}:=\sum_{n=0}^{\infty}\vert a_n\vert\exp\left(-cn^{1/(q+1)}
\right)<\infty$, is a countably normed  Fr\'{e}chet algebra.

Notice that for each $p>1$, $F^p$ is the Fr\'{e}chet envelope of the 
Privalov space $N^p$. In this paper we study the structure of maximal 
ideals in the algebras $F^p$ ($1<p<\infty$). In particular,
we give a complete characterization of closed maximal 
ideals in $F^p$. Moreover, we characterize 
multiplicative linear functionals on $F^p$.
   \end{abstract}

\maketitle

{\renewcommand{\thefootnote}{}

\vspace{-2mm}

\section{Introduction, Preliminaries and Results}

Let $\Bbb D$ denote the open unit disk 
in the complex plane  and let $\Bbb T$ denote the
boundary of $\Bbb D$.  Let $L^q(\Bbb T)$ $(0<q\le \infty)$ be the
familiar Lebesgue space on the unit circle $\Bbb T$.

The  Privalov class $N^p$ $(1<p<\infty)$ is defined 
as the set of all holomorphic functions $f$ $f$ on $\Bbb D$ such that 
   $$
\sup_{0<r<1}\int_0^{2\pi}(\log^+\vert f(re^{i\theta})\vert)^p\,
\frac{d\theta}{2\pi}<+\infty\eqno(1)
   $$ 
holds, where $\log^+|a|=\max\{\log |a|,0\}$.
 These classes were firstly considered  by Privalov in 
\cite[pages 93--10]{p}, where $N^p$ is denoted as $A_q$. 

Notice that for $p=1$, the condition (1) defines the 
 Nevanlinna class $N$ of holomorphic functions in $\Bbb D$.
 Recall that the  Smirnov class $N^+$  is the set
of all functions  $f$ holomorphic on $\Bbb D$ such that
   $$
\lim_{r\rightarrow 1}\int_0^{2\pi}\log^+\vert f(re^{i\theta})\vert
\,\frac{d\theta}{2\pi}=\int_0^{2\pi}\log^+\vert f^*(e^{i\theta})\vert\,
\frac{d\theta}{2\pi}<+\infty\eqno(2)
   $$
where $f^*$ is the boundary function of $f$ on $\Bbb T$; that is,
 $$
 f^*(e^{i\theta})=\lim_{r\rightarrow 1-}f(re^{i\theta})\eqno(3)
  $$ 
is the  radial  limit of $f$  which exists for almost 
every $e^{i\theta}\in\Bbb T$. 
 We denote by $H^q$ $(0<q\le\infty)$ the classical  Hardy space on 
$\Bbb D$. 

It is known (see \cite{mo, mp1, i1})  that the following inclusion 
relations hold: 
   $$
N^r\subset N^p\;(r>p),\quad\bigcup_{q>0}H^q\subset
\bigcap_{p>1}N^p,\quad{\rm and}\quad\bigcup_{p>1}N^p\subset N^+
\subset N,\eqno(4)
   $$
where the above containment relations are proper.
   
The study of the spaces $N^p$ $(1<p<\infty)$ was continued in 1977
by M. Stoll  \cite{s} (with the notation
 $(\log^+H)^\alpha$ in \cite{s}).  Further, the topological and functional 
properties of these spaces  have been  studied by several authors
(see \cite{mo}, \cite{e1},  \cite{e2},  \cite{im} and  
\cite{me4}--\cite{m5}).

M. Stoll \cite[Theorem 4.2]{s}  proved that for each  
 $p>1$ the space $N^p$ (with the noatation $(\log^+H)^{\alpha}$ 
in \cite{s}) equipped with the topology given by the metric 
$d_p$ defined by            
   $$
d_p(f,g)=\Big(\int_0^{2\pi}\big(\log(1+
\vert f^*(e^{i\theta})-g^*(e^{i\theta})\vert)\big)^p\,\frac{d\theta}
{2\pi}\Big)^{1/p},\quad f,g\in N^p,\eqno (5) 
   $$  
becomes an $F$-algebra, that is, $N^p$ is an  $F$-space (a complete 
metrizable topological vector space with the invariant metric)  in which  
multiplication is continuous.

 Recall that the function $d_1=d$ defined on the Smirnov class $N^+$
by (5) with $p=1$ induces the metric topology on $N^+$.
N. Yanagihara \cite{y1} showed that 
under this topology, $N^+$ is an $F$-space.

In connection with the spaces $N^p$
$(1<p<\infty)$,  Stoll \cite{s} 
(see also \cite{e1} and \cite[Section 3]{mp3})  also studied the spaces
$F^q$ $(0<q<\infty)$ (with the notation $F_{1/q}$ in \cite{s}), 
consisting of those functions $f$ holomorphic on 
$\Bbb D$ such that    
   $$
 \lim_{r\to 1}(1-r)^{1/q}\log^+M_{\infty}(r,f)=0,\eqno (6) 
   $$ 
where
$$
 M_{\infty}(r,f)=\max_{\vert z\vert\le r}\vert f(z)\vert\quad 
(0<r<1).\eqno (7) 
 $$

Here, as always in the sequel, we will need some Stoll's results
concerning the spaces 
$F^q$ only with $1<q<\infty$, and hence, we will assume that 
$q=p>1$ be any fixed number.

\vspace{2mm}

\noindent{\bf Theorem 1} (see \cite[Theorem 2.2]{s}).   {\it Suppose that
 $f(z)=\sum_{n=0}^\infty a_nz^n$ is a holomorphic function on $\Bbb D$. 
Then the following statements are equivalent:}
\begin{itemize}
\item[(a)] $f\in F^p$;
\item[(b)] {\it there exists a sequence $\{c_n\}_n$ of positive real
numbers with $c_n\to 0$ such that}
   $$
|a_n|\le \exp\left(c_nn^{1/(p+1)}\right),\quad n=0,1,2,\ldots;
\eqno(8) 
   $$
\item[(c)] {\it for any} $c>0$,
          $$
\Vert f\Vert_{p,c}:=\sum_{n=0}^{\infty}\vert a_n\vert\exp\left(-cn^{1/(p+1)}
\right)<\infty.\eqno(9) 
  $$   
\end{itemize}
\vspace{1mm}

\noindent{\it Remark 2.} Notice that in view of Theorem 1 
((a)$\Leftrightarrow$(c)), by (10) it is well defined  the family of seminorms 
$\left\{\Vert \cdot\Vert_{p,c} \right\}_{c>0}$  on $F^p$.

\vspace{2mm}

Recall that  a locally convex $F$-space is called 
a {\it Fr\'echet space}, and a {\it Fr\'echet algebra} 
is a  Fr\'echet space that is an algebra in which  multiplication
is continuous.  Stoll \cite{s} also proved the following result.

\vspace{2mm}

\noindent{\bf Theorem 3} (see \cite[Theorem 3.2]{s}).   {\it
The space $F^q$ $(0<q<\infty)$ equipped 
with the topology given by the  family of seminorms 
$\left\{\Vert \cdot\Vert_{q,c}\right\}_{c>0}$ defined for 
$f\in F^q$ as 
  $$
\Vert f\Vert_{q,c}:=\sum_{n=0}^{\infty}\vert a_n\vert\exp\left(-cn^{1/(q+1)}
\right)<\infty, \eqno (10) 
  $$
is a countably normed Fr\'{e}chet algebra}.

\vspace{2mm}

For our purposes, we will need the following result which 
characterizes the topological dual  of the space $F^p$.
   \vspace{2mm}

\noindent{\bf Theorem 4} (see \cite[Theorem 3.3]{s}).  {\it If $\gamma$ is a 
continuous linear 
functional on $F^p$, then there exists a sequence $\{\gamma_n\}_n$ of 
complex numbers with 
  $$
\gamma_n=O\big(\exp\left(-cn^{1/(p+1)}\right)\big),\quad
\mathrm{for\,\, some}\,\, c>0,\eqno(11)
   $$
such that 
               $$ 
\gamma(f)=\sum_{n=0}^{\infty}a_n\gamma_n,\eqno(12)
     $$        
 where $f(z)=\sum_{n=0}^{\infty}a_nz^n\in F^p$, with convergence
being absolute. Conversely, if $\{\gamma_n\}_n$ is a sequence of complex
 numbers for which 
 $$
\gamma_n=O\left(\exp\left(-cn^{1/(p+1)}\right)\right),\eqno(13)
$$
 then {\rm (13)} defines a continuous linear functional 
on $F^p $.}
\vspace{2mm}
 
 Notice  that the Privalov space $N^p$ $(1<p<\infty)$ is not locally 
convex (see \cite[Theorem 4.2]{e1} and \cite[Corollary]{m4}), 
and hence, $N^p$ is properly contained in $F^p$.
Moreover, $N^p$ is not locally bounded (see \cite[Theorem 1.1]{mp4}).
Moreover, Stoll showed  (\cite[Theorem 4.3]{s}) that for each $p>1$ 
$N^p$ is a dense subspace of $F^p$ and the topology on $F^p$ defined by the
family of seminorms (10)  is weaker than the topology on  $N^p$ 
given by the metric $d_p$ defined by} (5).
Furthermore,  Eoff  showed
\cite[Theorem 4.2, the case $p>1$]{e1} that $F^p$ is the
{\it Fr\'{e}chet envelope} of $N^p$. For more information on 
 Fr\'{e}chet envelope, see \cite[Theorem 1]{sh}, \cite[Section 1]{m2} and 
  \cite[Corollary 22.3, p. 210]{ke}.

 \vspace{2mm}
\noindent{\it Remark 5.}  
For $p=1$, the space $F_1$ has been denoted by $F^+$ and has been studied by 
N. Yanagihara in \cite{y2, y1}.  
It was shown in \cite{y2,y1} that $F^+$ is actually 
the containing Fr\'{e}chet space for $N^+$, 
i.e., $N^+$ with the initial topology embeds densely into $F^+$, under
the natural inclusion, and $F^+$ and the Smirnov class $N^+$ have the same 
topological duals.
\vspace{1mm}

Observe that the space $F^p$ topologised by the family of
seminorms $\left\{\Vert \cdot\Vert_{p,c}\right\}_{c>0}$ 
given by (10) is metrizable by the metric $\lambda_p$ defined as 
    $ \lambda_p(f,g)=\sum_{n=1}^{\infty}2^{-n}\frac{\Vert f-g\Vert_{p,1/n^
{p/(p+1)}}}{1+\Vert f-g\Vert_{p,1/n^{p/(p +1)}}}$ with $f,g\in F^p.
     $

Since Privalov space $N^p$ and its 
Fr\'{e}chet envelope $F^p$ $(1<p<\infty)$ are algebras, they can be 
also  considered as  rings with respect  to the usual  
ring's operations addition and multiplication. Notice that 
these two operations are continuous on $N^p$ and $F^p$  
because the spaces $N^p$ and $F^p$ become $F$-algebras. 

Motivated by several results on the ideal structure of some spaces 
of holomorphic functions given in \cite{k1} \cite{mo}, \cite{ma} and 
\cite{be}-\cite{Mor}, 
related investigations for the spaces $N^p$ $(1<p<\infty)$
and their Fr\'{e}chet envelopes were given in 
\cite{mo},  \cite{me4}, \cite{ma}, \cite{mp5}, \cite{mp3} and \cite{m5}.
Note that a survey of these results was given in \cite{mp6}.
The  $N^p$-analogue of the famous  Beurling's theorem 
for the Hardy spaces $H^q$ $(0<q<\infty)$ \cite{be} was proved in 
\cite{mp5}.
Moreover, it was proved in \cite[Theorem B]{me4})
that $N^p$ $(1<p<\infty)$ is a ring of Nevanlinna--Smirnov type  
in the sense of Mortini. The structure of closed weakly 
dense ideals in $N^p$ was established in \cite{mp3}.
The ideal structure of $N^p$ and the  multiplicative linear functionals  on 
$N^p$ were studied in \cite{mo} and \cite[Theorem ]{m5}. These results are
 similar to those obtained by Roberts and Stoll \cite{RS} for the 
Smirnov class $N^{+}$. 

Motivated by results of Roberts and Stoll given in  \cite[Section 2]{RS2} 
concerning a characterization of  multiplicative linear 
functionals on $F^{+}$ and closed maximal ideals in $F^{+}$,
in this paper we prove  the analogous results for 
the spaces $F^p$ $(1<p<\infty)$ given by Proposition 5, Proposition 6, 
Theorem 7 and Theorem 8.

\vspace{2mm}

\noindent{\bf Proposition 5.} {\it Let $\lambda\in\Bbb D$ and let 
$\gamma_\lambda$ be a functional on $F^p$ defined as
   $$
\gamma_\lambda(f)=f(\lambda)\eqno(14)
   $$ 
for every $f\in F^p$.  Then $\gamma_\lambda$ is a continuous multiplicative 
linear functional on $F^p$.}
\vspace{2mm}

For $\lambda\in\Bbb D$, we define 
  $$
{\mathcal M}_{\lambda}=\{ f\in F^p: f(\lambda)=0\}.\eqno(15)
  $$ 

\vspace{2mm}

\noindent{\bf Proposition 6.} 
{\it The set ${\mathcal M}_{\lambda}$ defined by $(15)$
is a closed maximal ideal in $F^p$ for each $\lambda\in \Bbb D$.}

\vspace{2mm}

\noindent{\bf Theorem 7.} {\it Let $\gamma$ be a nontrivial multiplicative
linear functional on $F^p$. Then there exists $\lambda\in \Bbb D$ such that
$$
\gamma(f)=f(\lambda)\eqno(16)
$$ 
for every $f\in F^p$. Moreover, $\gamma$ is 
a continuous map.}

\vspace{2mm}

\noindent{\bf Theorem 8.} {\it Let $p>1$ and let  
${\mathcal M}$ be a closed maximal ideal in $F^p$. Then there exists 
$\lambda\in \Bbb D$ such that ${\mathcal M}={\mathcal M}_{\lambda}$.}
\vspace{2mm}

\section{Proof of the results}

\vspace{2mm}

\noindent{\it Proof of Proposition 5.} Clearly, for each 
$\lambda\in\Bbb D$ $\gamma_{\lambda}$ is a multiplicative linear functional
on the space $F^p$. In order to show that $\gamma_{\lambda}$ is 
a continuous functional on $F^p$,  note that for 
 any function $f(z)=\sum_{n=0}^{\infty}a_nz^n\in F^p$ ($z\in\Bbb D$) we have 
  $$
\gamma_{\lambda}(f)=f(\lambda)=\sum_{n=0}^{\infty}a_n\lambda^n.\eqno(17)
  $$
Clearly,  the sequence $\{\gamma_n\}$ defined as 
$\gamma_n=\lambda^n$ ($n=0,1,2\ldots$) satisfies the asymptotic condition 
(11) of Theorem 4. This together with the equality (17) implies that 
$\gamma_{\lambda}$ is a continuous functional on $F^p$, and
the proof is completed. \hfill $\square$ 

\vspace{2mm}

\noindent{\it Proof of Proposition 6.}
Notice that in view of  (15), ${\mathcal M}_{\lambda}$
is the kernel of the functional $\gamma_{\lambda}$ 
defined on $F^p$ by (14). From this and the fact that by Proposition 5,
 $\gamma_{\lambda}$ is a continuous multiplicative linear functional
on the space $F^p$, we conclude that ${\mathcal M}_{\lambda}$ 
is a closed maximal ideal in $F^p$.  \hfill $\square$

\vspace{2mm}

\noindent{\it Proof of Theorem 7.} 
 If we take $\gamma(z)=\lambda$, then $\gamma(z-\lambda)=0$. 
If we suppose that $\lambda\notin \Bbb D$,
then $z\mapsto 1/(z-\lambda)$ ($z\in\Bbb D$) is a bounded 
function on the closed unit disk $\overline{\Bbb D}: |z|\le 1$.
Therefore, $z\mapsto z-\lambda$ ($z\in\Bbb D$) is an invertible element 
of the algebra $F^p$. If $f$ is any invertible element in $F^p$,
then $1=\gamma(1)=\gamma(f)\gamma(f^{-1})$, and thus, $\gamma(f)\not=0$.
Especially, we have $\gamma(z-\lambda)\not=0$. A contradiction, and 
hence, it must be $\lambda\in \Bbb D$. Then consider the set 
  $$
(z-\lambda)F^p:=\{(z-\lambda)f(z): f\in F^p\}.\eqno(18)
  $$  
For each $\lambda\in\Bbb D$, let ${\mathcal M}_\lambda$ be a set defined by 
(15).
Then obviously, $(z-\lambda)F^p\subset {\mathcal M}_\lambda$. 
 Moreover, if $f\in {\mathcal M}_\lambda$, then by (6) and (7) easily follows 
that $f$ can be expressed as a product $f(z)=(z-\lambda)g(z)$ with $g\in F^p$.
Therefore,
    $$
{\mathcal M}_\lambda=(z-\lambda)F^p,\eqno(19)
   $$ 
whence it follows that 
  $$
{\mathcal M}_\lambda\subseteq \ker\gamma ,\eqno(20)
   $$
where $\ker\gamma$ denotes the kernel of the functional $\gamma$. 
By Proposition 6, ${\mathcal M}_\lambda$ is a closed maximal ideal in 
$F^p$. This together with the  inclusion relation  (20) implies that 
  ${\mathcal M}_\lambda= \ker\gamma$. Moreover,  
$\gamma(f)=f(\lambda)$ for all $f\in F^p$ and $\gamma$ is continuous
on $F^p$ by Proposition 5. This completes the proof of the theorem. 
\hfill $\square$ 

 \vspace{2mm}

\noindent{\it Proof of Theorem 8.} We proceed as in 
\cite[Theorem 2]{St2}. If we set $X= F^p/{\mathcal M}$, 
then  in the terminology of Arens \cite{ar},  
$X$ is  complete, metrizable, convex complex topological 
division algebra. Therefore, by \cite{ar}, $X\cong \Bbb  C$. 
Thus, there exists a multiplicative linear functional 
$\gamma$ on  $F^p$ such that ${\mathcal M}=\ker\gamma$. 
Then by Theorem 7,
${\mathcal M}={\mathcal M}_{\lambda}$ for some  
$\lambda\in \Bbb D$, as asserted.
\hfill $\square$

\end{document}